\documentclass[11pt]{article}
\usepackage{amssymb, amsmath,  amsfonts,dsfont, mathrsfs,euscript,eufrak}
\usepackage{textcomp}
\usepackage{latexsym}
\usepackage{indentfirst}
\usepackage[normalem]{ulem}

\usepackage{graphicx}
\usepackage[colorlinks=true, allcolors=blue]{hyperref}

\newtheorem{Theorem}{Theorem}
\newtheorem{Lemma}{Lemma}
\newtheorem{Proposition}{Proposition}
\newtheorem{Corollary}{Corollary}

\newcommand{\dd}{\,d}

\newcommand{\R}{\mathbb R}
\newcommand{\Z}{\mathbb Z}
\newcommand{\N}{\mathbb N}

\newcommand{\RID}{\boldsymbol{Q}}

\newcommand{\ID}{\boldsymbol{I}}

\newcommand{\e}{\varepsilon}

\newcommand{\sF}{\widetilde{F}}

\newcommand*\wbar[1]{%
	\hbox{%
		\vbox{%
			\hrule height 0.5pt 
			\kern0.3ex
			\hbox{%
				\kern-0.18em
				\ensuremath{#1}%
				\kern0em
			}%
		}%
	}%
}

\begin{document}
\title{On necessary conditions of rational-infinite divisibility for distributions with non-zero discrete parts}

\author{A. A. Khartov$^{1,2,3,}$\footnote{Email addresses: \texttt{khartov.a@iitp.ru}, \texttt{alexeykhartov@gmail.com}} }

\footnotetext[1]{Institute for Information Transmission Problems (Kharkevich Institute) of Russian Academy of Sciences, Bolshoy Karetny per. 19, build.1, 127051 Moscow, Russia.}
\footnotetext[2]{Smolensk State University, 4 Przhevalsky st., 214000 Smolensk, Russia. }
\footnotetext[3]{ITMO University, 49 Kronverksky Pr., 197101 Saint-Petersburg, Russia.}

\maketitle
\begin{abstract}
	We consider the new class $\boldsymbol{Q}$ of rational-infinitely (or quasi-infinitely) divisible distribution functions on the real line. By definition, $F\in \boldsymbol{Q}$ if there are  some infinitely divisible distribution functions $F_1$ and $F_2$ such that $F_1=F*F_2$, where ``$*$'' is the convolution. The characteristic function of such $F$ admits the L\'evy--Khintchine-type representation
	with a ``signed spectral measure''. The class  $\boldsymbol{Q}$ is a significant extension of the family of  infinitely divisible distribution functions and it have already found some applications in several areas. So there is an active interest in this class. In particular, a lot of  results have recently appeared on the problem of belonging to the class $\boldsymbol{Q}$ in terms of characteristic functions. In the paper, we continue this series of results by proposing two necessary conditions for distribution functions from $\boldsymbol{Q}$ with non-zero discrete parts. Namely, the characteristic functions of such a distribution function and its discrete part are always separated from zero.	
\end{abstract}


\section{Introduction}

Let $F$ be a  distribution function on the real line $\R$. Let $f$ be the characteristic function of $F$, i.e.
\begin{eqnarray*}
	f(t):=\int_{\R} e^{itx} \dd F(x),\quad t\in\R.
\end{eqnarray*}
We call a distribution function $F$ (and the corresponding probability law) \textit{rational-infinitely divisible} if there are some infinitely divisible distribution functions $F_1$ and $F_2$ on $\R$ such that $F_1=F*F_2$, where ``$*$'' denotes the convolution operation. This equality may be written in the form $f(t)=f_1(t)/f_2(t)$, $t\in\mathbb{R}$, where $f_1$ and $f_2$ denote the characteristic functions of $F_1$ and $F_2$, respectively. Let us denote by $\RID$ the class of all rational-infinitely divisible distribution functions. It is clear that $\RID$ is an extension of the class $\ID$ of all infinitely divisible distribution functions. There are some examples of distribution functions from $\RID\setminus \ID$ in the classical monographs   \cite{GnedKolm} (p.~81--83) and  \cite{LinOstr} (p.~165).

The classes $\ID$ and $\RID$ have very similar properties. So it is seen from the definition that the characteristic function $f$ of any $F\in\RID$ doesn't have zeroes on $\R$ and it admits \textit{the L\'evy--Khinchine-type representation}:
\begin{eqnarray}\label{repr_f}
	f(t)=\exp\biggl\{it \gamma+\int_{\mathbb{R}} \Bigl(e^{itx} -1 -it \sin(x)\Bigr)\dfrac{1+x^2 }{x^2} \,d G(x)\biggr\},\quad t\in\mathbb{R},
\end{eqnarray}
with the \textit{shift parameter} $\gamma\in\mathbb{R}$, and  the \textit{spectral function} $G:\R\to\R$ of bounded total variation on $\R$. The function $G$ is assumed to be right-continuous at every point of the real line with the condition $G(-\infty)=0$.  Note that $G$ is non-monototic in general. Nevertheless, the \textit{spectral pair }$(\gamma, G)$ is uniquely determined by $f$ and hence by $F$ as for the class $\ID$ (see \cite{GnedKolm}, p.~80). Next, let $f$ be a characteristic function represented by formula \eqref{repr_f} with some $\gamma$ and $G$ satisfying the above-mentioned conditions. Following Lindner and Sato \cite{LindSato}, the corresponding distribution function $F$ (and the corresponding probability law) is called \textit{quasi-infinitely divisible}. Due to the Hahn--Jordan decomposition for $G$, it is not difficult to show that $F$ is rational-infinitely divisible.

The class $\RID$ is actively studied now. The main reference is the paper \cite{LindPanSato} by Lindner, Pan, and Sato, which contains the detailed analysis of this class. Multivariate analog of $\RID$ is considered in \cite{AlexeevKhartov_Mult}, \cite{BergKutLind},  and \cite{BergerLindner}. There are interesting applications in theory of stochastic processes (see  \cite{LindSato} and \cite{Pass}), number theory (see \cite{Nakamura}),  physics (see \cite{ChhDemniMou} and \cite{Demni}), and insurance mathematics (see  \cite{ZhangLiuLi}). The class $\RID$ penetrates statistics (see \cite{Panov}). 

This short note is connected with the problem of belonging to the class $\RID$. The general statement of the problem  is to obtain the necessary and/or sufficient conditions for which a given distribution function $F$ belongs to $\RID$. For instance, there is the following simple fact (see \cite{Cuppens}, Theorem~4.3.7): if a probability law has a mass $>1/2$ at some point, then its distribution function belongs to $\RID$. By now a number of deeper results have accumulated. We will give a brief survey of them now. Let $F$ be a distribution function on $\R$. According to the Lebesgue decomposition theorem,  $F$ admits the representation:
\begin{eqnarray}\label{repr_F_Lebesgue}
	F(x) = c_d F_{d}(x) + c_a F_{a}(x)+ c_s F_s(x), \quad x\in \R, 
\end{eqnarray}
where $F_d$, $F_a$, and $F_s$ are discrete, absolutely continuous and continuous singular distribution functions, respectively. Here the coefficients $c_a$, $c_a$, and $c_s$  are non-negative constants such that $c_d+c_a+c_s=1$. The characteristic function  $f$ of $F$ is  represented accordingly:
\begin{eqnarray}\label{repr_f_Lebesgue}
	f(t) = c_d f_{d}(t) + c_a f_{a}(t)+ c_s f_s(t), \quad t\in \R, 
\end{eqnarray}
where $f_d$, $f_a$, and $f_s$ are the characteristic functions corresponding to $F_d$, $F_a$, and  $F_s$, respectively. The most of existing results are formulated in terms of characteristic functions. The first rather general result was obtained by Lindner, Pan, and Sato in \cite{LindPanSato} (Theorem~8.1, p.~30) for the case $c_d=1$ and discrete lattice $F=F_d$. It states that such $F$ belongs to the class $\RID$  if and only if  $f(t)\ne 0$ for any $t\in\R$. The criterion for arbitrary discrete $F$ (exactly the case $c_d=1$) was obtained in \cite{KhartovSPL2022} (see also \cite{AlexeevKhartov}): $F\in\RID$ if and only if $f$ is separeted from zero, i.e. $\inf_{t\in\R}|f(t)|>0$.  This criterion generalizes  the previous one, because the absolute value $|f(\,\cdot\,)|$  is a periodic continuous function on $\R$ for lattice $F$ and therefore, in the lattice case, $f$ is zero-free on the period segment (and hence on $\R$) if any only if it is separated from zero.
By the way, there is a recent result in paper \cite{KhartovKharchenkov}, where this criterion is formulated in terms of the support points of the discrete distribution and their probabilities (see also \cite{KhartovAlexeevThree}). Next, in the paper \cite{BergerKutlu}, Berger and Kutlu considered the case $c_d>0$, $c_a\geqslant 0$, and   $c_s=0$.  They proved that, under these assumptions, $F\in\RID$ if and only if \textbf{1)} $f(t)\ne 0$ for any $t\in\R$ and \textbf{2)} $\inf_{t\in\R}|f_d(t)|>0$. Conditions  \textbf{1)} and  \textbf{2)} together are equivalent to one condition \textbf{3)} $\inf_{t\in\R}|f(t)|>0$. We note that this result is a generalization of an earlier criterion by Berger (see \cite{Berger}), where $F_d$ is assumed to be lattice.  Next, it is asserted in \cite{KhartovMixSingPart} that the previous criterion is valid under a more general assumption of the dominated continuous singular part. More precisely, let us assume that $c_d>0$, $c_a\geqslant0$, $c_s\geqslant 0$, and define $\mu_d:=\inf_{t\in\R}|f_d(t)|$. We additionally assume that $c_s< c_d \mu_d$ if $\mu_d>0$ and $c_s=0$ if $\mu_d=0$. Then $F\in\RID$ if and only if conditions \textbf{1)} and \textbf{2)} hold together or equivalent condition  \textbf{3)} is satified. It is important to note, however, that there is an interesting example of $F\notin \RID$, for which $c_s=\mu_d c_d>0$, \textbf{1)} and \textbf{2)} are satisfied, but \textbf{3)} is not (see \cite{KhartovMixSingPart}, Example~3). In other words,   the pair of  \textbf{1)} and \textbf{2)} is not equivalent to \textbf{3)} in general case and both  $\inf_{t\in\R}|f_d(t)|>0$ and $\inf_{t\in\R}|f(t)|>0$ are important for belonging to the class $\RID$ in the case $c_d>0$. In addition to those listed, there are several results where other approaches are used (see \cite{BergerLindner}, \cite{KhartovLogarithm}, and \cite{KhartovGeneral}).

In this note,  we propose two results, where we will deal with the case $c_d>0$ too. We will show that, actually, the conditions  $\inf_{t\in\R}|f_d(t)|>0$ and $\inf_{t\in\R}|f(t)|>0$ (without any additional assumptions) are always  necessary for belonging $F$ to the class $\RID$. We hope that these results will make some progress to the most general potential criterion of membership of the class $\RID$ in the case $c_d>0$ or they will be some useful tools for further results of other kind in this theme.

\section{Results}

Let $F$ be a (right-continuous) distribution function on $\R$ with the characteristic function $f$ as before. For further convenience we rewrite formulas \eqref{repr_F_Lebesgue} and \eqref{repr_f_Lebesgue} in the following forms
\begin{eqnarray}\label{repr_Ffc}
	F(x)=c_dF_d(x)+(1-c_d) F_c(x),\quad x\in\R,\quad\text{and}\quad f(t)=c_d f_d(t)+(1-c_d) f_c(t),\quad t\in\R,
\end{eqnarray}
respectively. Here $c_d\in[0,1]$, $F_c$ denotes a continuous distribution function and $f_c$ denotes its characteristic function. If $c_d=1$ then $F_c$ can be chosen arbitrarily. If $c_d<1$, i.e. $c_a+c_s=1-c_d>0$, then
\begin{eqnarray*}
	F_c(x)=\dfrac{c_a F_a(x)+c_s F_s(x)}{c_a+c_s},\quad x\in\R,\quad\text{and}\quad f_c(t)=\dfrac{c_a f_a(t)+c_s f_s(t)}{c_a+c_s},\quad t\in\R.
\end{eqnarray*}
Note that we don't assume a priori that $F\in\RID$ (namely, in Lemma~\ref{lm_MeanCont}, Proposition~\ref{pr_1} and \ref{pr_2} below). 

Let us formulate the first main result of the paper.

\begin{Theorem}\label{th_fsep0}
	If $F\in\RID$ and $c_d>0$, then $\inf_{t\in\R}|f(t)|>0$.
\end{Theorem}

We will need the following simple lemmata to prove this and further assertions.

\begin{Lemma}\label{lm_fraction_tpmh}
	Suppose that  $F\in\RID$. Then  there are constants $B\geqslant 0$ and $C>0$ such that
	\begin{eqnarray*}
		\dfrac{|f(t-h)f(t+h)|}{|f(t)|^2}\leqslant C e^{B h^2}\quad \text{for any}\quad t\in\R,\,\,  h\in\R.
	\end{eqnarray*}	
\end{Lemma}
\textbf{Proof of Lemma~\ref{lm_fraction_tpmh}.} Since $F\in\RID$, the function $f$ admits representation \eqref{repr_f} with some $\gamma$ and $G$. Using this formula, it is not difficult to show that
\begin{eqnarray*}
	\dfrac{f(t-h)f(t+h)}{f(t)^2}=\exp\biggl\{-\int_{\mathbb{R}}e^{itx} \bigl(1-\cos(hx)\bigr)\dfrac{1+x^2 }{x^2} \,d G(x)  \biggr\}
\end{eqnarray*}	
for any real $t$ and $h$ (see \cite{KhartovWeak}, p.~359, or \cite{KhartovGeneral}, p.~4, for more details). 
 \begin{eqnarray*}
 	\dfrac{|f(t-h)f(t+h)|}{|f(t)|^2}&=&\exp\biggl\{-\int_{\mathbb{R}}\cos(tx) \bigl(1-\cos(hx)\bigr)\dfrac{1+x^2 }{x^2} \dd G(x)  \biggr\}\\
 	&\leqslant& \exp\biggl\{\int_{\mathbb{R}}\bigl(1-\cos(hx)\bigr)\dfrac{1+x^2 }{x^2} \dd|G|(x)  \biggr\},
 \end{eqnarray*}
where $|G|(x)$ denotes the total variation of $G$ on $(-\infty,x]$ for any $x\in\R$. It is known that $(1-\cos(hx))\tfrac{1+x^2}{x^2}\leqslant \tfrac{h^2}{2}+2$ for any real $x$ and $h$ (see \cite{KhartovWeak}, p.~358). Therefore
 \begin{eqnarray*}
	\dfrac{|f(t-h)f(t+h)|}{|f(t)|^2}\leqslant \exp\biggl\{\biggl(\dfrac{h^2}{2}+2\biggr)\cdot \|G\| \biggr\}\quad \text{for any}\quad t\in\R,\,\,  h\in\R,
\end{eqnarray*}
where $\|G\|$ denotes the total variation of $G$ on $\R$. Defining $B:=\tfrac{1}{2}\,\|G\|\geqslant 0$ and $C:=e^{2\|G\|}>0$, we come to the required estimate. \quad $\Box$

\begin{Lemma}\label{lm_MeanCont}
	It is always valid that
	\begin{eqnarray*}
		\sup_{t\in\R} \biggl\{\dfrac{1}{2T}\int_{-T}^{T}|f_c(t+h)|^2 \dd h\biggr\}\to 0\quad\text{as} \quad T\to\infty.
	\end{eqnarray*}		
\end{Lemma}
\textbf{Proof of Lemma~\ref{lm_MeanCont}.} By definition, $f_c$ is the characteristic function for $F_c$. Hence the function $t\mapsto |f_c(t)|^2$, $t\in\R$, is the characteristic function of the symmetrized version of $F_c$, say $\sF_c$ (see \cite{Lukacs}, p.~38). Thus we have
\begin{eqnarray*}
	|f_c(t)|^2=\int_{\R} e^{itx} \dd \sF_c(x),\quad t\in\R.
\end{eqnarray*}
Let us consider the following family of integrals:
\begin{eqnarray*}
	M(t,T):=\dfrac{1}{2T}\int_{-T}^{T}|f_c(t+h)|^2\dd h=\dfrac{1}{2T}\int_{-T}^{T}\biggl(\int_{\R} e^{i(t+h)x} \dd \sF_c(x)\biggr)\dd h,\quad t\in\R,\,\, T>0.
\end{eqnarray*}
It is clear that we can write
\begin{eqnarray*}
	M(t,T)=\int_{\R} e^{itx}\biggl(\dfrac{1}{2T}\int_{-T}^{T} e^{ihx} \dd h\biggr)\dd \sF_c(x),\quad t\in\R,\,\, T>0.
\end{eqnarray*}
Then
\begin{eqnarray*}
	M(t,T)&=&\int_{\{0\}}\biggl(\dfrac{1}{2T}\int_{-T}^{T}\dd h\biggr) \dd \sF_c(x)  + \int_{\R\setminus\{0\}}e^{itx}\biggl(\dfrac{1}{2T}\int_{-T}^{T} e^{ihx} \dd h\biggr)\dd \sF_c(x)\\
	&=&\int_{\{0\}} \dd \sF_c(x)+\int_{\R\setminus\{0\}}e^{itx}\,\dfrac{\sin(Tx)}{Tx}\dd \sF_c(x),\quad t\in\R,\,\, T>0.
\end{eqnarray*}
Here $\int_{\{0\}} \dd \sF_c(x)= \sF_c(0)-\sF_c(0-)=0$, because $\sF_c$, being the simmetrization of the continuous distribution function $F_c$, is continuous too. Consequently, for any $T>0$
\begin{eqnarray*}
	\sup_{t\in\R}|M(t,T)|= \sup_{t\in\R}\biggl| \int_{\R\setminus\{0\}}e^{itx}\,\dfrac{\sin(Tx)}{Tx}\dd \sF_c(x)\biggr|\leqslant  \int_{\R\setminus\{0\}}\biggl|\dfrac{\sin(Tx)}{Tx}\biggr|\dd \sF_c(x).
\end{eqnarray*}
Observe that, in the latter integral, the integrand function is bounded by $1$ and it vanishes as $T\to\infty$. Therefore this integral vanishes too  by Lebesgue's dominated convergence theorem. Thus we obtain $\sup_{t\in\R}|M(t,T)|\to 0$ as $T\to\infty$ as required.\quad $\Box$\\

Now we are ready to prove our first result.\\

\noindent\textbf{Proof of Theorem~\ref{th_fsep0}.}  By assumption, $F\in\RID$ and hence we know that $f(t)\ne 0$ for any $t\in\R$ (see introduction).  In order to prove $\inf_{t\in\R}|f(t)|>0$, we will use the quotient $f(t-h)f(t+h)/f(t)^2$ with real $t$ and $h$. Namely, let us consider the following integral
\begin{eqnarray}\label{def_I}
	I(t,\tau):=\dfrac{1}{2\tau}\int_{-\tau}^{\tau}\dfrac{|f(t-h)f(t+h)|}{|f(t)|^2}\,\dd h
\end{eqnarray}
for any $t\in\R$ and $\tau>0$. According to Lemma~\ref{lm_fraction_tpmh}, this integral admits the estimate
\begin{eqnarray*}
	I(t,\tau)\leqslant \dfrac{1}{2\tau}\int_{-\tau}^{\tau} C e^{B h^2}\dd h\leqslant C e^{B\tau^2}<\infty
\end{eqnarray*}
with some constants $B\geqslant 0$ and $C>0$. In particular, we see that
\begin{eqnarray}\label{conc_supIttau}
	\sup_{t\in\R} I(t,\tau)<\infty\quad\text{for any}\quad \tau>0.
 \end{eqnarray}

We now assume, contrary to our claim, that $\inf_{t\in\R}|f(t)|=0$. We will show below that for some $\tau>0$ it yields $\sup_{t\in\R} I(t,\tau)=\infty$,  which contradicts \eqref{conc_supIttau}. For this purpose we will obtain a lower estimate for the integral
\begin{eqnarray}\label{def_J}
	J(t,\tau):= I(t,\tau)\cdot|f(t)|^2=\dfrac{1}{2\tau}\int_{-\tau}^{\tau}|f(t-h)f(t+h)|\,\dd h
\end{eqnarray}
for any $t\in\R$ and $\tau>0$. 

Let us consider the integrand function in \eqref{def_J}. Using the decomposition for $f$ from \eqref{repr_Ffc},  for any $t\in\R$ and $h\in\R$ we write
\begin{eqnarray*}
	f(t-h)f(t+h)&=&\bigl(c_d f_d(t-h)+(1-c_d) f_c(t-h) \bigr)\cdot \bigl(c_d f_d(t+h)+(1-c_d) f_c(t+h) \bigr)\\
	&=&c_d^2 f_d(t-h)f_d(t+h)+(1-c_d)^2 f_c(t-h)f_c(t+h)\\
	&&+c_d(1-c_d) f_d(t+h)f_c(t-h)+c_d(1-c_d) f_d(t-h)f_c(t+h).
\end{eqnarray*}
So it follows that
\begin{eqnarray*}
	|f(t-h)f(t+h)|&\geqslant&c_d^2 |f_d(t-h)f_d(t+h)|-(1-c_d)^2 |f_c(t-h)f_c(t+h)|\\
	&&-c_d(1-c_d) |f_d(t+h)f_c(t-h)|-c_d(1-c_d)| f_d(t-h)f_c(t+h)|\\
	&\geqslant& c_d^2 |f_d(t-h)f_d(t+h)|- (1-c_d)^2|f_c(t-h)f_c(t+h)|\\
	&&{}-c_d(1-c_d) |f_c(t-h)|-c_d(1-c_d)|f_c(t+h)|,
  \end{eqnarray*}
where we used the inequality $|f_d(s)|\leqslant 1$, $s\in\R$. Then for any $t\in\R$ and $\tau>0$ 
 \begin{eqnarray*}
 	J(t,\tau)&\geqslant&  \dfrac{c_d^2}{2\tau}\int_{-\tau}^{\tau}|f_d(t-h)f_d(t+h)|\,\dd h- \dfrac{(1-c_d)^2}{2\tau}\int_{-\tau}^{\tau}|f_c(t-h)f_c(t+h)|\,\dd h\\
 	&&- \dfrac{c_d(1-c_d)}{2\tau}\int_{-\tau}^{\tau}|f_c(t-h)|\,\dd h-\dfrac{c_d(1-c_d)}{2\tau}\int_{-\tau}^{\tau}|f_c(t+h)|\,\dd h.
 \end{eqnarray*}
It is easily seen that
 \begin{eqnarray*}
 	\dfrac{1}{2\tau}\int_{-\tau}^{\tau}|f_c(t-h)|\,\dd h=\dfrac{1}{2\tau}\int_{-\tau}^{\tau}|f_c(t+h)|\,\dd h.
 \end{eqnarray*}
Due to the inequality $|f_c(s)|\leqslant 1$ for any $s\in\R$, we get
 \begin{eqnarray*}
 	\dfrac{1}{2\tau}\int_{-\tau}^{\tau}|f_c(t-h)f_c(t+h)|\,\dd h\leqslant  \dfrac{1}{2\tau}\int_{-\tau}^{\tau}|f_c(t+ h)|\,\dd h.
 \end{eqnarray*} 
Therefore
\begin{eqnarray*}
	J(t,\tau)\geqslant \dfrac{c_d^2}{2\tau}\int_{-\tau}^{\tau}|f_d(t-h)f_d(t+h)|\,\dd h-\dfrac{(1-c_d)^2}{2\tau}\int_{-\tau}^{\tau}|f_c(t+h)|\,\dd  h-\dfrac{c_d(1-c_d)}{\tau}\int_{-\tau}^{\tau}|f_c(t+h)|\,\dd  h.
\end{eqnarray*}
Thus for any $t\in\R$ and $\tau>0$ we obtain 
 \begin{eqnarray}\label{ineq_JJdJc}
 	J(t,\tau)\geqslant\dfrac{c_d^2}{2\tau}\int_{-\tau}^{\tau}|f_d(t-h)f_d(t+h)|\,\dd h-\dfrac{1-c_d^2}{2\tau}\int_{-\tau}^{\tau}|f_c(t+h)|\,\dd  h\geqslant c_d^2 J_d(t,\tau)- J_c(t,\tau),
 \end{eqnarray}
where we define
\begin{eqnarray*}
	J_d(t,\tau):=\dfrac{1}{2\tau}\int_{-\tau}^{\tau}|f_d(t-h)f_d(t+h)|\,\dd h,\qquad J_c(t,\tau):= \dfrac{1}{2\tau}\int_{-\tau}^{\tau}|f_c(t+h)|\,\dd  h,\quad t\in\R,\,\,\tau>0.
\end{eqnarray*}

Next, according to the assumption $\inf_{t\in\R}|f(t)|=0$, we choose a sequence $(t_k)_{k\in\N}$ such that $|f(t_k)|\to  0$ as $k\to\infty$ (here and below, $\N$ denotes the set of positive integers). Since $f(t)\ne 0$ for any $t\in\R$ and $f$ is a continuous function on $\R$, it follows that $|t_k|\to \infty$ as $k\to\infty$. Without loss of generality,
we can  assume that $t_k\to \infty$ as $k\to\infty$, because $|f(t)|=|f(-t)|$, $t\in\R$. Next, using the sequence $(t_k)_{k\in\N}$, we define the functions $\varphi_k(h):= f_d(t_k+h)$, $h\in\R$, $k\in\N$. Since $f_d$ is an almost periodic function, the sequence $(\varphi_k)_{k\in\N}$ is relatively compact in the topology of uniform convergence on $\R$ (see \cite{Levitan}, pp. 23--24). So there exists a subsequence $(\varphi_{k_l})_{l\in\N}$ that  uniformly converges to an almost periodic function $\varphi$, i.e.
\begin{eqnarray}\label{conc_conv}
	\sup_{h\in\R}|f_d(t_{k_l}+h)-\varphi(h)|\to 0,\quad l\to\infty.
\end{eqnarray}
From this we have 
\begin{eqnarray}\label{conc_conv_pmphi}
	\sup_{h\in\R}|f_d(t_{k_l}+h)f_d(t_{k_l}-h)-\varphi(h)\varphi(-h)|\to 0,\quad l\to\infty.
\end{eqnarray}
Here the function $h\mapsto\varphi(h)\varphi(-h)$, $h\in\R$, is an almost periodic (as a product of two almost periodic functions,  see \cite{Levitan}, p.~27) and  it is known that there exists $h_0\in\R$ such that $\varphi(h_0)\varphi(-h_0)\ne 0$ (see \cite{KhartovSPL2022}, p.~3, for more details).  In addition, by the Parseval identity for almost periodic functions (see \cite{KhartovSPL2022}, p.~45), the following limit exists 
\begin{eqnarray}\label{def_C}
	A:=\lim\limits_{T\to\infty} \dfrac{1}{2T} \int_{-T}^{T} |\varphi(h)\varphi(-h)|^2 \dd h,
\end{eqnarray}   
and $A$ is the sum of  squared absolute values of  all the Fourier coefficients of  $\varphi(h)\varphi(-h)$, $h\in\R$. Since this function is not identically zero, we conclude  that $A>0$ (see \cite{KhartovSPL2022}, p.~45).

We now turn to the sequence
\begin{eqnarray*}
	J_d(t_{k_l},\tau)=\dfrac{1}{2\tau}\int_{-\tau}^{\tau}|f_d(t_{k_l}-h)f_d(t_{k_l}+h)|\,\dd h,\quad l\in\N,
\end{eqnarray*}
with arbitarary $\tau>0$. The following inequalities hold
 \begin{eqnarray*}
 	J_d(t_{k_l},\tau)&\geqslant& \dfrac{1}{2\tau}\int_{-\tau}^{\tau}|\varphi(h)\varphi(-h)|\,\dd h-\dfrac{1}{2\tau}\int_{-\tau}^{\tau}|f_d(t_{k_l}+h)f_d(t_{k_l}-h)-\varphi(h)\varphi(-h)|\,\dd h\\
 	&\geqslant&\dfrac{1}{2\tau}\int_{-\tau}^{\tau}|\varphi(h)\varphi(-h)|\,\dd h- \sup_{h\in\R}\bigl|f_d(t_{k_l}+h)f_d(t_{k_l}-h)-\varphi(h)\varphi(-h)\bigr|.
 \end{eqnarray*}
Let us fix arbitrary $\e\in(0,\tfrac{1}{4})$. Due to \eqref{conc_conv_pmphi}, there exists $l_\e\in\N$ such that for any integer $l\geqslant l_\e$
\begin{eqnarray*}
	\sup_{h\in\R} |f_d(t_{k_l}+h)f_d(t_{k_l}-h)-\varphi(h)\varphi(-h)|<\e A,
\end{eqnarray*}
and, consequently, we have 
\begin{eqnarray*}
	J_d(t_{k_l},\tau)\geqslant\dfrac{1}{2\tau}\int_{-\tau}^{\tau}|\varphi(h)\varphi(-h)|\,\dd h- \e A.
\end{eqnarray*}
Next, we know that $|f_d(s)|\leqslant 1$ for any $s\in\R$. Due to \eqref{conc_conv}, the same inequality is true for $\varphi$, i.e. $|\varphi(s)|\leqslant 1$ for any $s\in\R$. Then $|\varphi(\pm h)|\geqslant |\varphi(\pm h)|^2$ for any $h\in\R$ and we have
 \begin{eqnarray}\label{ineq_Jd}
	J_d(t_{k_l},\tau)\geqslant\dfrac{1}{2\tau}\int_{-\tau}^{\tau}|\varphi(h)\varphi(-h)|^2\,\dd h- \e A
\end{eqnarray}
for any integer $l\geqslant l_\e$ and $\tau>0$. 

We next consider the integrals $J_c(t_{k_l},\tau)$, $l\in\N$, $\tau>0$. It is
 convenient to proceed with the squared integrand:
\begin{eqnarray*}
	J_c(t_{k_l},\tau)= \dfrac{1}{2\tau}\int_{-\tau}^{\tau}|f_c(t_{k_l}+h)|\,\dd  h\leqslant \sqrt{ \dfrac{1}{2\tau}\int_{-\tau}^{\tau}|f_c(t_{k_l}+h)|^2\,\dd  h },\quad l\in\N,\,\, \tau>0.
\end{eqnarray*} 
 
We now choose a suitable value of $\tau$. According to \eqref{def_C}, for the fixed $\e$ there exists $\tau_{d,\e}>0$ such that for any $\tau\geqslant\tau_{d,\e}$
\begin{eqnarray*}
	\dfrac{1}{2\tau}\int_{-\tau}^{\tau}|\varphi(h)\varphi(-h)|^2\,\dd h\geqslant (1-\e)A. 
\end{eqnarray*}
Hence, on account \eqref{ineq_Jd}, we have $J_d(t_{k_l},\tau)\geqslant (1-\e) A-\e A=(1-2\e)A$ for any integer $l\geqslant l_\e$ and $\tau\geqslant \tau_{d,\e}$. Next, by Lemma~\ref{lm_MeanCont},   for the fixed $\e$ we find $\tau_{c,\e}>0$ such that for any $\tau\geqslant\tau_{c,\e}$
\begin{eqnarray*}
	\sup_{t\in\R}\biggl\{\dfrac{1}{2\tau}\int_{-\tau}^{\tau}|f_c(t+h)|^2\,\dd h\biggr\}\leqslant c_d^4\,\e^2 A^2. 
\end{eqnarray*}
It is important  that here $c_d>0$ by the assumption of the theorem. Consequently, $J_c(t_{k_l},\tau)\leqslant c_d^2\,\e A$ for any $l\in\N$ and $\tau\geqslant\tau_{c,\e}$. Let $\tau_\e:=\max\{\tau_{d,\e},\tau_{c,\e}\}$. Then $J_d(t_{k_l},\tau_{\e})\geqslant (1-2\e)A$  and $J_c(t_{k_l},\tau_{\e})\leqslant c_d^2\,\e A$ are satisfied for any integer $l\geqslant l_\e$. We combine these inequalities in \eqref{ineq_JJdJc}:
\begin{eqnarray*}
	J(t_{k_l},\tau_{\e})\geqslant c_d^2 J_d(t_{k_l},\tau_{\e})-J_c(t_{k_l},\tau_{\e})\geqslant c_d^2(1-2\e) A-c_d^2\,\e A= c_d^2(1-3\e) A.
\end{eqnarray*}
Recall that $\e<\tfrac{1}{4}$. So, in particular, $J(t_{k_l},\tau_{\e})>\tfrac{1}{4} c_d^2 A>0$.

Let us return to the function $I$ defined by \eqref{def_I}. According to \eqref{def_J} and the above estimate, we have
\begin{eqnarray*}
	I(t_{k_l},\tau_{\e})=\dfrac{J(t_{k_l},\tau_{\e})}{|f(t_{k_l})|^2}>\dfrac{c_d^2 A}{4|f(t_{k_l})|^2}
\end{eqnarray*}
for any integer $l\geqslant l_\e$. There is the strictly positive constant $c_d^2 A$ in the numerator, but $|f(t_{k_l})|\to 0$ as $l\to \infty$ by the assumption. Hence $	I(t_{k_l},\tau_{\e})\to \infty$ as $l\to\infty$, which contradicts \eqref{conc_supIttau} as we expected. Thus $\inf_{t\in\R}|f(t)|=0$ can not be valid and we come to the opposite, i.e. $\inf_{t\in\R}|f(t)|>0$.\quad$\Box$

\begin{Proposition}\label{pr_1}
	If $\inf_{t\in\R}|f(t)|>0$ then $c_d>0$.
\end{Proposition}
\textbf{Proof of Proposition~\ref{pr_1}.} Contrary to our claim, suppose that $c_d=0$. Hence $F=F_c$ and $f=f_c$ according to \eqref{repr_Ffc}. Then, on the one hand, 
\begin{eqnarray*}
	\dfrac{1}{2T}\int_{-T}^{T} |f(h)|^2\dd h\to 0\quad\text{as}\quad T\to\infty
\end{eqnarray*}
by Lemma~\ref{lm_MeanCont}. On the other hand, by the assumption of the theorem,  $\mu_d:=\inf_{t\in\R}|f(t)|>0$, and hence
\begin{eqnarray*}
	\dfrac{1}{2T}\int_{-T}^{T} |f(h)|^2\dd h\geqslant 	\dfrac{1}{2T}\int_{-T}^{T} \mu_d^2\dd h=\mu_d^2\quad\text{for any}\quad T>0.
\end{eqnarray*}
It contradicts to the previous convergence and thus $c_d=0$ is not possible, i.e. we have $c_d>0$. \quad $\Box$\\

On account of this proposition, we can modify the statement of Theorem~\ref{th_fsep0} in the following way.

\begin{Corollary}
	Suppose that $F\in\RID$. Then $\inf_{t\in\R} |f(t)|>0$ if and only if $c_d>0$.
\end{Corollary}

The next assertion significantly complement Proposition~\ref{pr_1}.

\begin{Proposition}\label{pr_2}
	If $\inf_{t\in\R}|f(t)|>0$ then $\inf_{t\in\R}|f_d(t)|>0$.
\end{Proposition}
\textbf{Proof of Proposition~\ref{pr_2}.} If $c_d=1$, i.e. $f=f_d$, the assertion is obvious. We focus on the case $c_d\in(0,1)$. Let $\mu:=\inf_{t\in\R}|f(t)|>0$. Suppose, contrary to our claim, that
\begin{eqnarray}\label{assum_fd}
	\inf_{t\in\R}|f_d(t)|=0.
\end{eqnarray}
So we fix arbitrary $\e \in\bigl(0,\tfrac{1}{4}\bigr)$ and we choose $t_{\e}\in\R$ such that $|f_d(t_\e)|<\e\mu$. Since $f_d$ is an almost periodic function,  there exists an $\ell_\e>0$  such that every segment $[a,a+\ell_\e]$ with $a\in\R$ contains at least one number $\tau_{\e}$, for which the inequality $|f_d(t+\tau_{\e})-f_d(t)|<\e\mu$ holds for any $t\in\R$ (see \cite{Levitan}, p.~20). Let $\Z$ denote the set of integers as usual. So we introduce the following family of segments
\begin{eqnarray*}
	S_{\e,k}:=\bigl[\bigl(2k-\tfrac{1}{2}\bigr)\ell_\e, \bigl(2k+\tfrac{1}{2}\bigr)\ell_\e\bigr],\quad k\in\Z,
\end{eqnarray*}
and for every $k\in\Z$ we choose a point $\tau_{\e,k}\in S_{\e,k}$ as above, i.e.
\begin{eqnarray}\label{ineq_supfdttau}
	\sup_{t\in\R}\bigl|f_d(t+\tau_{\e,k})-f_d(t)\bigr|<\e\mu,\quad k\in\Z.
\end{eqnarray} 
Since $0\in S_{\e,0}$, it is clear that we  can set $\tau_{\e,0}:=0$. Next, by the continuity of $f_d$ on $\R$, we find $\delta_\e\in \bigl(0,\tfrac{1}{2}\ell_\e\bigr]$ such that $|f_d(t_\e+\delta)-f_d(t_\e)|<\e\mu$ for any $\delta\in[0,\delta_\e)$. Due to \eqref{ineq_supfdttau}, we have
\begin{eqnarray*}
	\bigl|f_d(t_\e+\delta+\tau_{\e,k})-f_d(t_\e+\delta)\bigr|<\e\mu,\quad k\in\Z,
\end{eqnarray*} 
for any $\delta\in[0,\delta_\e)$. Next, we define the sequence $t_{\e,k}:=t_\e+\tau_{\e,k}$, $k\in\Z$, where $t_{\e,0}=t_\e$. According to the definitions, the numbers $t_{\e,k}$ are distinct and, more precisely, 
\begin{eqnarray*}
	\ldots< t_{\e,-2}<t_{\e,-1}<t_{\e,0}<t_{\e,1}<t_{\e,2}<\ldots
\end{eqnarray*}
Moreover,  the intervals $(t_{\e,k}-\delta_\e,t_{\e,k}+\delta_\e)$, $k\in\Z$, are disjoint, i.e.
\begin{eqnarray}\label{ineq_tek}
	t_{\e,k}+\delta_\e<t_{\e,k+1}-\delta_\e,\quad k\in\Z.
\end{eqnarray}
Let us estimate the absolute values of $f_d$ on these intervals. We first observe that for any $\delta\in[0,\delta_\e)$ and $k\in\Z$
\begin{eqnarray*}
	\bigl|f_d(t_{\e,k}+\delta)\bigr|&\leqslant&  \bigl|f_d(t_\e+\delta)\bigr|+\bigl|f_d(t_{\e,k}+\delta)-f_d(t_\e+\delta)\bigr|\\
	&\leqslant& \bigl|f_d(t_\e)\bigr|+ \bigl|f_d(t_\e+\delta)-f_d(t_\e)\bigr|+\bigl|f_d(t_{\e,k}+\delta)-f_d(t_\e+\delta)\bigr|.
\end{eqnarray*}
Hence, by the above inequalities, we have   $|f_d(t_{\e,k}+\delta)|<3\e\mu$ for any $\delta\in[0,\delta_\e)$ and $k\in\Z$.

We now turn to the function $f_c$ from the decomposition. Due to $c_d<1$,  it is expressed by the formula $f_c(t)=\tfrac{f(t)-c_df_d(t)}{1-c_d}$, $t\in\R$. So we observe that
\begin{eqnarray*}
	|f_c(t_{\e,k}+\delta)|=\biggl|\dfrac{f(t_{\e,k}+\delta)-c_df_d(t_{\e,k}+\delta)}{1-c_d}\biggr|\geqslant \dfrac{|f(t_{\e,k}+\delta)|-c_d|f_d(t_{\e,k}+\delta)|}{1-c_d}.
\end{eqnarray*}
Hence the obtained estimates yield
\begin{eqnarray}\label{ineq_fctekdelta}
	|f_c(t_{\e,k}+\delta)|\geqslant \dfrac{\mu-3c_d\e\mu}{1-c_d}>(1-3\e)\mu
\end{eqnarray}
for any $\delta\in[0,\delta_\e)$ and $k\in\Z$. We next consider the following sequence of integrals
\begin{eqnarray*}
	I_{\e,n}:=\dfrac{1}{2T_{\e,n}}\int_{-T_{\e,n}}^{T_{\e,n}}|f_c(t_\e+h)|^2 \dd h,\quad n\in\N_0,
\end{eqnarray*}	
with $T_{\e,n}:=(2n+1)\ell_\e$, $n\in\N_0$, where $\N_0:=\N\cup \{0\}$. On the one hand, since $T_{\e,n}\to \infty$ as $n\to\infty$, we know that
\begin{eqnarray}\label{conv_Ien}
	I_{\e,n}\to 0\quad\text{as} \quad n\to\infty
\end{eqnarray}
by Lemma~\ref{lm_MeanCont}. On the other hand,  for any $n\in\N_0$ we have 
\begin{eqnarray}\label{conc_UnTen}
	U_{\e,n}:=\bigcup_{k=-n}^n (t_{\e,k}-\delta_\e,t_{\e,k}+\delta_\e)\subset \bigl[t_\e- T_{\e,n}, t_\e+T_{\e,n}\bigr].
\end{eqnarray}
Indeed, it follows from $\eqref{ineq_tek}$  and the inequalities 
\begin{eqnarray*}
	t_{\e,n}+\delta_\e=t_\e+\tau_{\e,n}+\delta_\e< t_\e+\bigl(2n+\tfrac{1}{2}\bigr) \ell_e+\tfrac{1}{2}\ell_\e=t_\e+(2n+1)\ell_\e=t_\e+T_{\e,n},\\
	t_{\e,-n}-\delta_\e=t_\e+\tau_{\e,-n}-\delta_\e> t_\e+\bigl(-2n-\tfrac{1}{2}\bigr) \ell_e-\tfrac{1}{2}\ell_\e=t_\e-(2n+1)\ell_\e=t_\e-T_{\e,n}.
\end{eqnarray*}
Due to \eqref{conc_UnTen}, we have
\begin{eqnarray*}
	I_{\e,n}=\dfrac{1}{2T_{\e,n}}\int_{t_\e-T_{\e,n}}^{t_\e+T_{\e,n}}|f_c(h)|^2 \dd h\geqslant \dfrac{1}{2T_{\e,n}}\int_{U_{\e,n}}|f_c(h)|^2 \dd h=\dfrac{1}{2T_{\e,n}}\sum_{k=-n}^{n}\int_{t_{\e,k}-\delta_\e}^{t_{\e,k}+\delta_\e}|f_c(h)|^2 \dd h,\quad n\in\N_0,
\end{eqnarray*}
where, on account of \eqref{ineq_fctekdelta},
\begin{eqnarray*}
	\int_{t_{\e,k}-\delta_\e}^{t_{\e,k}+\delta_\e}|f_c(h)|^2 \dd h=\int_{-\delta_\e}^{\delta_\e}|f_c(t_{\e,k}+\delta)|^2 \dd \delta\geqslant (1-3\e)^2\mu^2 \int_{-\delta_\e}^{\delta_\e}\dd \delta=2\delta_\e (1-3\e)^2\mu^2.
\end{eqnarray*}
Then
\begin{eqnarray*}
	I_{\e,n}\geqslant \dfrac{1}{2T_{\e,n}} \sum_{k=-n}^{n} 2\delta_\e (1-3\e)^2\mu^2 = \dfrac{(2n+1)\delta_\e (1-3\e)^2\mu^2}{T_{\e,n}} =\dfrac{\delta_\e (1-3\e)^2\mu^2}{\ell_\e}>\dfrac{\delta_\e \mu^2}{16\ell_\e},
\end{eqnarray*}
where we used that $\e\in\bigl(0,\tfrac{1}{4}\bigr)$. So we see that $I_{\e,n}\nrightarrow 0$ as $n\to\infty$ and we come to a contradiction with \eqref{conv_Ien}. Thus \eqref{assum_fd} is impossible and we conclude that $\inf_{t\in\R}|f_d(t)|>0$.\quad $\Box$\\
                                                                           
Theorem~\ref{th_fsep0} and Proposition~\ref{pr_2} yield at once the following assertion, which is the second main result of the paper.
  
\begin{Theorem}\label{th_fdsep0}
	If $F\in\RID$ and $c_d>0$, then $\inf_{t\in\R}|f_d(t)|>0$.
\end{Theorem}

Thus we proved that the conditions $\inf_{t\in\R}|f(t)|>0$ and $\inf_{t\in\R}|f_d(t)|>0$ are necessary for belonging $F$ to the class $\RID$ in the case $c_d>0$.  Since the first condition is always stronger than the second one (as we showed above) and they are not equivalent in general case (see introduction), only  the first inequality should be included in a  group of potential sufficient conditions for generalization of the mentioned existing results.

\end{document}